\documentstyle{article}

\newcommand{\bigzerou}{%
\smash{\lower1.7ex\hbox{\bg 0}}}
\newtheorem{defi}{Definition}

\newtheorem{conj}{Conjecture}

\newcommand{\ba}{\begin{eqnarray}}
\newcommand{\ea}{\end{eqnarray}}
\newcommand{\no}{\nonumber}

\newcommand{\mapright}[1]{%
\smash{\mathop{%
\hbox to 1.0cm{\rightarrowfill}}\limits^{#1}}}
\newcommand{\mapleft}[1]{%
\smash{\mathop{%
\hbox to 1.3cm{\leftarrowfill}}\limits^{#1}}}

\begin{document}
\title{
\begin{flushright}
  \begin{minipage}[b]{5em}
    \normalsize
    ${}$      \\
  \end{minipage}
\end{flushright}
{\bf Completion of the Conjecture: Quantum Cohomology of Fano Hypersurfaces}}
\author{Masao Jinzenji\\
\\ 
\it Graduate School of Mathematical Sciences\\
\it University of Tokyo\\
\it  Meguro-ku, Tokyo 153-8914, Japan\\}
\maketitle
\begin{abstract}
In this paper, we propose the formulas that compute all the
rational structural constants of the quantum K\"ahler sub-ring 
of Fano hypersurfaces.  
\footnote{ e-mail 
address: jin@hep-th.phys.s.u-tokyo.ac.jp}  
\end{abstract}
\section{Introduction and Statement of the Main Results}
The main result of this paper is the completion of the conjectures 
proposed in \cite{cj} on the quantum K\"ahler sub-ring of Fano hypersurfaces. 

The quantum cohomology ring of Fano hypersurfaces are studied both in 
physical and mathematical points of view \cite{mp}. 
One of the most successful approach from
mathematics is the use 
of Gauss-Manin system with the aid of the gravitational descendants 
\cite{giv}, \cite{kim}, \cite{bert}. 
These works showed that the quantum cohomology 
ring of Fano hypersurfaces are realized as the Gauss-Manin connection 
of the Gauss-Manin system, whose solution is given by the generating 
function of a certain type of the gravitational correlation functions 
(including the gravitational descendants) of the topological sigma model 
coupled to the topological gravity. Moreover, they also showed that 
the solutions of the Gauss-Manin system are written as the
generalized hyper geometric functions in the case of 
Fano hypersurfaces. Their result is beautiful, but in compensation
for the introduction of gravitational descendants, it is not very clear to 
see the explicit structural constants and the multiplication 
table of the quantum cohomology ring from their results.

Here, we construct a procedure that computes the explicit
structural constants and the multiplication table of the quantum
cohomology ring of Fano hypersurfaces without using  
the gravitational descendants.

In \cite{cj}, we have constructed the recursive formulas 
that describe the rational structural constant $L_{n}^{N,k,d}$ of the 
quantum K\"ahler sub-ring of 
$M_{N}^{k}$, which is the degree $k$ hypersurfaces in $CP^{N-1}$, in 
terms of those of $M_{N+1}^{k}$ when the degree of the rational curves 
concerned is no more than $5$. These recursive formulas 
determine the structural 
constants of the quantum cohomology ring of $M_{N}^{k}$ with the aid of 
the initial conditions given by the following formula in \cite{ej}:     
\begin{equation}
\sum_{n=0}^{k-1}L_{n}^{N,k,1}w^{n}=k\prod_{j=1}^{k-1}(jw+(k-j)),
\label{one}
\end{equation}
where $N-k$ must be greater than $1$.
But in \cite{cj}, our construction of the recursive formula
needs some explicit data of the structural constants coming from 
the higher degree curves, and the 
determination of the recursive formulas gets 
extremely harder as the degree of the rational curves grows.  

In this paper, we propose the formulas that give the recursive 
formula for the rational curves of arbitrary degree $d$
without using any numerical 
data of the explicit structural constants. 
Combining this procedure with the formula (\ref{one}), 
we can compute all the structural constants of the 
quantum K\"ahler sub-ring of $M_{N}^{k}\; (N\geq k)$
by use of our previous results in \cite{cj}.

Here, we summarize the result to be proposed in this paper. 
First, we introduce the formal variable 
$x,y,z_{1},z_{2},\cdots,z_{d-1},t_{1},t_{2},\cdots,t_{d-1}$ ($d$ is the 
degree of the rational curves in the hypersurface). 
Then we take the following linear combination of the variables 
$x,y,t_{1},t_{2},\cdots,t_{d-1}$:  
 \begin{equation}
\varphi_{i}(x,y,t_{1},t_{2},\cdots,t_{d-1})=\frac{(d-i)\cdot x+i\cdot y}{d}+
\sum_{j=1}^{i-1}\frac{d-i}{d-j}\cdot t_{j}+
\sum_{j=i+1}^{d-1}\frac{i}{j}\cdot t_{j}.
\end{equation}
Next, consider the rational function $p_{i}$ in the formal variables,
\begin{equation}
p_{i}(x,y,z_{i},t_{1},t_{2},\cdots,t_{d-1})=t_{i}+\varphi_{i}
+z_{i}\cdot\frac{t_{i}+\varphi_{i}}{t_{i}+\varphi_{i}-z_{i}}.
\label{p}
\end{equation}
And we expand $\frac{t_{i}+\varphi_{i}}{t_{i}+\varphi_{i}-z_{i}}$ 
into the form,
\begin{equation}
\sum_{n=0}^{\infty}\sum_{k=0}^{\infty}{-n\choose k}
(\frac{z_{i}}{t_{i}})^{n}(\frac{\varphi_{i}}{t_{i}})^{k},
\label{expand}
\end{equation}
where
\begin{eqnarray} 
{-n\choose 0}:=1,\quad
{-n\choose k}:=\frac{1}{k!}\prod_{j=1}^{k}(-n-j+1)\quad \;(k\geq1).
\end{eqnarray}
Then we take the product of these formal series
$p_{i}=t_{i}+\varphi_{i}+
z_{i}(\sum_{n=0}^{\infty}\sum_{k=0}^{\infty}{-n\choose k}
(\frac{z_{i}}{t_{i}})^{n}(\frac{\varphi_{i}}{t_{i}})^{k})$,
\begin{equation}
R_{d}(x,y,z_{1},z_{2},\cdots,z_{d-1},t_{1},t_{2},\cdots,t_{d-1})=
\prod_{j=1}^{d-1}p_{j}(x,y,z_{j},t_{1},t_{2},\cdots,t_{d-1}),
\label{fact}
\end{equation}
and pick up the terms that are free from (constant terms in) $t_{i}$'s,
\begin{eqnarray}
&&Poly_{d}(x,y,z_{1},z_{2},\cdots,z_{d-1})\no\\
&&=\frac{1}{(2\pi\sqrt{-1})^{d-1}}
\int_{C_{1}}\frac{dt_{1}}{t_{1}}\cdots
\int_{C_{d-1}}\frac{dt_{d-1}}{t_{d-1}}R_{d}(x,y,z_{1},z_{2},\cdots,z_{d-1},
t_{1},\cdots ,t_{d-1}).
\label{residue}
\end{eqnarray}
In (\ref{residue}), we have to choose the path $C_{i}$ in the domain
where the expansion in (\ref{expand}) is effective.
 $Poly_{d}$ is the polynomial first introduced in the formula
(5.96) in \cite{cj} (web-site version), and it is possible to reconstruct 
the recursive formula for the degree $d$ rational curves from it.
 We prepare some definitions 
that are needed for the reconstruction of 
the recursive formula from $Poly_{d}$.
Consider the monomial 
$x^{d_{i_{0}}}z_{i_{1}}^{d_{i_{1}}}\cdots z_{i_{m}}^{d_{i_{m}}}
y^{d_{i_{m+1}}}
\;(\sum_{j=0}^{m+1}d_{i_{j}}=d-1)$, that appear in $Poly_{d}$,
associated with the following ``comb type'' of a positive integer $d$ 
\cite{bert}:
\begin{equation}
0=i_{0}<i_{1}<i_{2}<\cdots<i_{m}<i_{m+1}=d.
\end{equation}
Next, we prepare some elements in (a free 
abelian group) ${\bf Z}^{m+1}$, which are 
determined for each monomial 
$x^{d_{i_{0}}}z_{i_{1}}^{d_{i_{1}}}\cdots z_{i_{m}}^{d_{i_{m}}}
y^{d_{i_{m+1}}}$ ,
as follows:
\begin{eqnarray}
\alpha&:=&(m+1-d,m+1-d,\cdots,m+1-d),\no\\
\beta&:=&(0,i_{1}-1,i_{2}-2,\cdots,i_{m}-m),\no\\
\gamma&:=&(0,i_{1}(N-k),i_{2}(N-k),\cdots,i_{m}(N-k)),\no\\
\epsilon_{1}&:=&(1,0,0,0,\cdots,0),\no\\
\epsilon_{2}&:=&(1,1,0,0,\cdots,0),\no\\
\epsilon_{3}&:=&(1,1,1,0,\cdots,0),\no\\
&&\cdots\no\\
\epsilon_{m+1}&:=&(1,1,1,1,\cdots,1).
\end{eqnarray}
Now we define $\delta=(\delta_{1},\cdots,\delta_{m+1})\in{\bf Z}^{m+1}$ 
by the formula: 
\begin{equation}
\delta:=\alpha+\beta+\gamma+\sum_{j=1}^{m}(d_{i_{j}}-1)\epsilon_{j}+
d_{i_{m+1}}\epsilon_{m+1}.
\label{delta}
\end{equation}
Then our main conjecture is the following:
\begin{conj}
The recursive formula for the degree $d$ rational curves is obtained from
$Poly_{d}$ by the formula:
\begin{equation}
L_{n}^{N,k,d}=\phi(Poly_{d}),
\end{equation}
where $\phi$ is a ${\bf Q}$-linear map from the ${\bf Q}$-vector
space of the homogeneous polynomials of degree $d-1$ in $x,y,z_{1},\cdots
,z_{d-1}$ to the ${\bf Q}$-vector
space of the weighted homogeneous polynomials of degree $d$ in 
$L_{m}^{N+1,k,d'}$ defined as follows: 
\begin{equation}
\phi(x^{d_{0}}y^{d_{d}}z_{i_{1}}^{d_{i_{1}}}\cdots z_{i_{m}}^{d_{i_{m}}})
=\prod_{j=1}^{m+1}L_{n+\delta_{j}}^{N+1,k,i_{j}-i_{j-1}}.
\end{equation}
\end{conj}

This paper is organized as follows.
In Section 2, we introduce the notation on the quantum K\"ahler sub-ring
of the hypersurface in $CP^{N-1}$ and briefly review the results
obtained in \cite{cj}.
In Section 3, we first observe
the agreement between the results obtained by 
Conjecture 1 and the recursive formulas up to 
the $d=5$ rational curves derived in \cite{cj}. Next, we briefly 
explain the trend of thoughts, that leads us to the formula proposed in 
Conjecture 1.
In Section 4, we test our conjecture in the $d=6$ case. We explicitly 
write down the recursive formula for the $d=6$ rational curves, which 
is derived from the direct use of Conjecture 1 and see that this formula 
has the non-trivial property conjectured in \cite{cj}. 
In Appendix A,
we present the table of the recursive formulas derived in \cite{cj} 
, that are heavily referred in the body of this paper. 

\section{Quantum K\"ahler Sub-ring of Projective Hypersurfaces}
\subsection{Notation} 
In this section, we introduce the quantum K\"ahler sub-ring 
of the quantum cohomology ring of a degree $k$ hypersurface in
$CP^{N-1}$.
Let $M_{N}^{k}$ be a hypersurface of degree $k$ in $CP^{N-1}$.
 We denote by $QH^{*}_{e}(M_{N}^{k})$ the 
subring of the quantum cohomology ring $QH^{*}(M_{N}^{k})$
generated by ${\cal O}_{e}$ induced from the K\"ahler form $e$ 
(or, equivalently the intersection $H\cap M_{N}^{k}$ between a hyperplane
class $H$ of $CP^{N-1}$ and $M_{N}^{k}$).
 The multiplication rule of $QH^{*}_{e}(M_{N}^{k})$ 
is determined by the Gromov-Witten invariant of genus $0$ 
$\langle {\cal O}_{e}{\cal O}_{{e}^{N-2-m}}
{\cal O}_{{e}^{m-1-(k-N)d}}\rangle_{d,M_{N}^{k}}$ and
it is given as follows:
\begin{eqnarray}
 L_{m}^{N,k,d} &:=&\frac{1}{k}\langle {\cal O}_{e}{\cal O}_{{e}^{N-2-m}}
{\cal O}_{{e}^{m-1-(k-N)d}}\rangle,\no\\
\no\\
{\cal O}_{e}\cdot 1&=&{\cal O}_{e},\nonumber\\
{\cal O}_{e}\cdot{\cal O}_{{e}^{N-2-m}}&=&{\cal O}_{{e}^{N-1-m}}+
\sum_{d=1}^{\infty}L_{m}^{N,k,d}q^{d}{\cal O}_{{e}^{N-1-m+(k-N)d}},\no\\
q&:=&\exp(t),
\label{gm}
\end{eqnarray}
where the subscript $d$ counts the degree of the rational curves
measured by $e$. So $q=\exp(t)$ is the degree counting 
parameter. 
Since $M_{N}^{k}$ is a complex $(N-2)$ dimensional manifold, we see that
a structure constant $L_{m}^{N,k,d}$
is non-zero only if the following condition is satisfied:
\begin{eqnarray}
&& 1\leq N-2-m\leq N-2, 1\leq m-1+(N-k)d\leq N-2,\no\\
&\Longleftrightarrow &max\{0,2-(N-k)d\}\leq m \leq min\{N-3,N-1-(N-k)d\}.
\label{sel}
\end{eqnarray}
We rewrite (\ref{sel}) into 
\begin{eqnarray}
(N-k\geq 2) &\Longrightarrow& 0\leq m \leq (N-1)-(N-k)d\no\\
(N-k=1,d=1)&\Longrightarrow& 1\leq m \leq N-3\no\\
(N-k=1,d\geq2)&\Longrightarrow& 0\leq m \leq N-1-(N-k)d\no\\
(N-k\leq 0)&\Longrightarrow& 2+(k-N)d\leq m \leq N-3.
\label{flasel}
\end{eqnarray}
From (\ref{flasel}), we easily see that the number of the non-zero
structure constants $L_{m}^{N,k,d}$ is finite except for the case of $N=k$.
Moreover, if $N\geq 2k$, the non-zero structure constants come only from
the $d=1$ part and the non-vanishing $L_{m}^{N,k,1}$  
is determined by $k$ and  
independent of $N$. 
The $N\geq 2k$ region is studied by Beauville \cite{beauville}, 
and his result plays 
the role of an initial condition of our discussion later.
Explicitly, they are given by the formula (\ref{one}).
In the case of $N=k$, the multiplication rule of $QH^{*}_{e}(M_{k}^{k})$ is
given as follows:
\begin{eqnarray}
{\cal O}_{e}\cdot 1&=&{\cal O}_{e},\nonumber\\
{\cal O}_{e}\cdot{\cal O}_{{e}^{k-2-m}}&=&
(1+\sum_{d=1}^{\infty}q^{d}L_{m}^{k,k,d}){\cal O}_{{e}^{k-1-m}}
\;\;(m=2,3,\cdots,k-3),\no\\
{\cal O}_{e}\cdot{\cal O}_{{e}^{k-3}}  &=&{\cal O}_{e^{k-2}}.
\label{calabi}
\end{eqnarray}
We introduce here the generating function of the structure constants 
of the Calabi-Yau hypersurface $M_{k}^{k}$:
\begin{equation}
L_{m}^{k,k}(e^{t}):=1+\sum_{d=1}^{\infty}L_{m}^{k,k,d}e^{dt}\;\;(m=2,\cdots 
,k-3).
\end{equation}
\subsection{Review of Results for Fano and Calabi-Yau Hypersurfaces
and Virtual Structure Constants}
Let us summarize the results of \cite{cj}. In \cite{cj}, 
we showed that the structure constants $L_{m}^{N,k,d}$ 
of $QH_{e}^{*}(M_{N}^{k})$ for $(N-k\geq 2)$  can be obtained by 
applying the recursive formulas in Appendix A (up to the $d=5$ case), 
which describe
$L_{m}^{N,k,d}$ in terms of $L_{m'}^{N+1,k,d'}\;\;(d'\leq d)$,
with the initial 
conditions $L_{m}^{N,k,1}$ 
and $L_{m}^{N,k,d}=0\;\;(d\geq 2)$ in the $N\geq 2k$ region.
These recursive formulas naturally lead us to the relation:
\begin{equation}
({\cal O}_{e})^{N-1}-k^{k}({\cal O}_{e})^{k-1}q=0
\label{jincolli}
\end{equation}
of $QH_{e}^{*}(M_{N}^{k})\:\:(N-k\geq 2)$            
by descending induction using Beauville's result \cite{beauville}, \cite{jin}. 

In the $N-k=1$ case, the recursive formulas receive modification 
only in the $d=1$ part:
\begin{equation}
L_{m}^{k+1,k,1}=L_{m}^{k+2,k,1}-L_{0}^{k+2,k,1}=L_{m}^{k+2,k,1}-k!.
\label{givgiv}
\end{equation}
This leads us to the following relation of $QH_{e}^{*}(M_{k+1}^{k})$:
\begin{equation}
({\cal O}_{e}+k!q)^{N-1}-k^{k}({\cal O}_{e}+k!q)^{k-1}q=0.
\label{ch1}
\end{equation} 
 
The structural constant $L_{m}^{k,k,d}$ for a Calabi-Yau hypersurface 
does not obey the recursive formulas in appendix A.    
We introduce here the virtual structure constants $\tilde{L}_{m}^{N,k,d}$
as follows.
\begin{defi}
Let $\tilde{L}_{m}^{N,k,d}$ be the rational number obtained by applying
the recursion relations of Fano hypersurfaces in Appendix A
 for arbitrary $N$ and $k$ with the initial condition
 $L_{n}^{N,k,1}\;\;(N\geq 2k)$ and 
$L_{n}^{N,k,d}=0\;\;(d\geq 2,\;\;N\geq 2k)$. 
\end{defi}
{\bf Remark 1}
In the $N-k\geq2$ region, $\tilde{L}_{m}^{N,k,d}=L_{m}^{N,k,d}$.

We define the generating function of the virtual structure constants of 
the Calabi-Yau hypersurface $M_{k}^{k}$
as follows:
\begin{eqnarray}
\tilde{L}^{k,k}_{n}(e^{x})&:=&1+
\sum_{d=1}^{\infty}\tilde{L}^{k,k,d}_{n}e^{dx},\nonumber\\
&&(n=0,1,\cdots,k-1).
\label{str2}
\end{eqnarray}
In \cite{cj}, we conjectured that $\tilde{L}^{k,k}_{n}(e^{x})$ gives us the 
information of the B-model of the mirror manifold of $M_{k}^{k}$. 
More explicitly, we conjectured 
\begin{eqnarray}
&&\tilde{L}_{0}^{k,k}(e^{x})=
\sum_{d=0}^{\infty}\frac{(kd)!}{(d!)^{k}}e^{dx},\no\\
&&\tilde{L}_{1}^{k,k}(e^{x})=
\frac{dt(x)}{dx}:=\frac{d}{dx}
\biggl(x+(\sum_{d=1}^{\infty}\frac{(kd)!}{(d!)^{k}}\cdot 
(\sum_{i=1}^{d}\sum_{m=1}^{k-1}\frac{m}{i(ki-m)})e^{dx})/
(\sum_{d=0}^{\infty}\frac{(kd)!}{(d!)^{k}}e^{dx})\biggr)
\label{po}
\end{eqnarray}
where the r.h.s. of (\ref{po}) is derived from the solutions of 
the ODE for the period integral of the mirror manifold of $M_{k}^{k}$,
\begin{equation}
((\frac{d}{dx})^{k-1}-ke^{x}(k\frac{d}{dx}+1) 
(k\frac{d}{dx}+2)\cdots(k\frac{d}{dx}+k-1))w(x)=0,
\label{ode1}
\end{equation}
that was used in the computation based on the mirror symmetry, 
\cite{candelas2}, \cite{gmp}, \cite{mnmj}.
Of course, we can extend the conjecture (\ref{po}) to the 
general $\tilde{L}_{n}^{k,k}(e^{x})$ if we compare the 
$\tilde{L}_{n}^{k,k}(e^{x})$ with the B-model three point functions
in \cite{gmp}.
Hence we obtain the mirror map $t=t(x)$ without 
using the mirror conjecture:
\begin{equation}
t(x)=x+\int_{-\infty}^{x}dx'({\tilde{L}^{k,k}_{1}(e^{x'})}-1)
=x+\sum_{d=1}^{\infty}\frac{\tilde{L}^{k,k,d}_{1}}{d}e^{dx}.
\end{equation}
With the conjecture given by (\ref{po}), we can construct 
the mirror transformation that transforms 
the virtual structure constants of the  Calabi-Yau hypersurface into 
the real ones as follows:
\begin{equation}
L^{k,k}_{m}(e^{t})=\frac{\tilde{L}^{k,k}_{m}(e^{x(t)})}
{\tilde{L}^{k,k}_{1}(e^{x(t)})}.\quad(m=2,\cdots,k-3) 
\label{mt}
\end{equation}
Therefore we can compute all the structural constants of quantum K\"ahler sub-ring 
of $M_{N}^{k}\;(N\geq k)$ if we can determine the form of the 
recursive formulas for the rational curves of arbitrary degree 
as in Appendix A.   
\section{Agreement with the Known Results}
In the first part of this section, we show that our main conjecture 
reproduces the recursive 
formulas in the  $d\leq 5$ cases obtained in \cite{cj}.
First, we write down $R_{d}$ in (\ref{fact}) up to the 
$d=5$ case. Here, we omit 
some irrelevant terms that do not appear in $Poly_{d}$.
\begin{eqnarray}
R_1&=&1\no\\
R_2&=&(\frac{x+y}{2}+z_1)\no\\
R_3&=&(\frac{2x+y}{3}+z_1+\frac{(z_{1})^{2}}{t_1}+\frac{1}{2}t_{2})
(\frac{x+2y}{3}+z_2+\frac{(z_{2})^{2}}{t_2}+\frac{1}{2}t_{1})\no\\
R_4&=&(\frac{3x+y}{4}+z_1+\frac{(z_{1})^{2}}{t_1}(1-\frac{3x+y}{4t_{1}})+
\frac{(z_{1})^{3}}{(t_1)^{2}}+\frac{1}{2}t_{2}+\frac{1}{3}t_{3})\times
\nonumber\\
&&(\frac{2x+2y}{4}+z_2+\frac{(z_{2})^{2}}{t_2}(1-\frac{2x+2y}{4t_{2}})+
   \frac{(z_{2})^{3}}{(t_2)^{2}}+\frac{2}{3}t_{1}+\frac{2}{3}t_{3})\times
\nonumber\\
&&(\frac{x+3y}{4}+z_3+\frac{(z_{3})^{2}}{t_3}(1-\frac{x+3y}{4t_{3}})+
   \frac{(z_{3})^{3}}{(t_3)^{2}}+\frac{1}{3}t_{1}+\frac{1}{2}t_{2})\no\\  
R_5&=&\biggl(\frac{4x+y}{5}+\frac{1}{2}t_{2}+\frac{1}{3}t_{3}
+\frac{1}{4}t_{4}+z_1+\no\\
&&\frac{(z_{1})^{2}}{t_1}(t_{1}/(t_{1}+\frac{4x+y}{5}+
\frac{1}{2}t_{2}+\frac{1}{3}t_{3}+\frac{1}{4}t_{4}))+
\frac{(z_{1})^{3}}{(t_1)^{2}}
(1-2(\frac{4x+y}{5t_{1}}))+\frac{(z_{1})^{4}}{(t_1)^{3}}\biggr)
\times\nonumber\\
&&\biggl(\frac{3x+2y}{5}+\frac{3}{4}t_{1}+\frac{2}{3}t_{3}+\frac{1}{2}t_{4}+
z_2+\nonumber\\
&&\frac{(z_{2})^{2}}{t_2}(t_{2}/(t_{2}+\frac{3x+2y}{5}+
\frac{3}{4}t_{1}+\frac{2}{3}t_{3}+\frac{1}{2}t_{4}))+
\frac{(z_{2})^{3}}{(t_2)^{2}}
(1-2(\frac{3x+2y}{5t_{2}}))+\frac{(z_{2})^{4}}{(t_2)^{3}}\biggr)
\times\nonumber\\
&&\biggl(\frac{2x+3y}{5}+\frac{1}{2}t_{1}+\frac{2}{3}t_{2}+\frac{3}{4}t_{4}+
z_3+\nonumber\\
&&\frac{(z_{3})^{2}}{t_3}(t_{3}/(t_{3}+\frac{2x+3y}{5}+
\frac{1}{2}t_{1}+\frac{2}{3}t_{2}+\frac{3}{4}t_{4}))+
\frac{(z_{3})^{3}}{(t_3)^{2}}
(1-2(\frac{2x+3y}{5t_{3}}))+\frac{(z_{3})^{4}}{(t_3)^{3}}\biggr)
\times\nonumber\\
&&\biggl(\frac{x+4y}{5}+\frac{1}{4}t_{1}+\frac{1}{3}t_{2}+\frac{1}{2}t_{3}+
z_4+\nonumber\\
&&\frac{(z_{4})^{2}}{t_4}(t_{4}/(t_{4}+\frac{x+4y}{5}+
\frac{1}{4}t_{1}+\frac{1}{3}t_{2}+\frac{1}{2}t_{3}))+
\frac{(z_{4})^{3}}{(t_4)^{2}}
(1-2(\frac{x+4y}{5t_{4}}))+\frac{(z_{4})^{4}}{(t_4)^{3}}\biggr)
\label{r}
\end{eqnarray}
In the $d\leq 4$ cases, the derivation of $Poly_{d}$ is just picking up 
the constant terms in ${t_{i}}$ from the above polynomials.
\begin{eqnarray}
Poly_1&=&1\no\\
Poly_2&=&(\frac{x+y}{2}+z_1)\no\\
Poly_{3}&=&(\frac{2x^{2}+5xy+2y^{2}}{9})
          +(\frac{2x+y}{3}+\frac{1}{2}z_{2})z_{2}
          +(\frac{x+2y}{3}+\frac{1}{2}z_{1})z_{1}
          +z_{1}z_{2}\nonumber\\
Poly_{4}&=&(\frac{3x^{3}+13x^{2}y+13xy^{2}+3y^{3}}{32})\nonumber\\
       &&+(\frac{x^{2}+4xy+3y^{2}}{8}
          +(\frac{3x+11y}{18})z_{1}
           +\frac{2}{9}(z_{1})^{2})z_{1}\nonumber\\
       &&+(\frac{3x^{2}+10xy+3y^{2}}{16}+(\frac{3x+3y}{8})z_{2} 
         +\frac{1}{4}(z_{2})^{2})z_{2}\nonumber\\
       &&+(\frac{3x^{2}+4xy+y^{2}}{8}+(\frac{11x+3y}{18})z_{3}
           +\frac{2}{9}(z_{3})^{2})z_{3}\nonumber\\
       &&+(\frac{3x+y}{4}+\frac{1}{2}z_{2}+\frac{1}{3}z_{3})z_{2}z_{3}
         +(\frac{x+y}{2}+\frac{2}{3}z_{1}+\frac{2}{3}z_{3})
         z_{1}z_{3}\nonumber\\
       &&+(\frac{x+3y}{4}+\frac{1}{3}z_{1}+\frac{1}{2}z_{2})
         z_{1}z_{2}\nonumber\\
       &&+z_{1}z_{2}z_{3}
\end{eqnarray}
Our main conjecture in Section 1 straightforwardly leads us 
to the recursive formulas in Appendix A.
But in the $d=5$ case, we have to take residues of some rational functions 
in $t_{i}$ in evaluating the coefficients of the monomials 
$(z_{i})^{2}(z_{j})^{2}\;\;(i<j)$. For example, we have to
evaluate the following residue integral to determine the coefficient of 
$(z_{1})^{2}(z_{4})^{2}$, 
\begin{equation}
\frac{1}{(2\pi\sqrt{-1})^{2}}
\int_{C_{1}}dt_{1}\int_{C_{4}}dt_{4}
\frac{(\frac{3}{4}t_{1}+\frac{1}{2}t_{4})
(\frac{1}{2}t_{1}+\frac{3}{4}t_{4})}
{t_{1}t_{4}(t_{1}+\frac{1}{4}t_{4})(t_{4}+\frac{1}{4}t_{1})}
=\frac{2}{3}.
\label{res}
\end{equation}
After these integrations, we can see the complete agreement
of $Poly_{d}$ derived from (\ref{r}) 
with the recursive formula for the $d=5$ rational curves.

Now, we briefly explain the process how we reached the formula 
$R_{d}$. Fundamental idea is the observation in \cite{cj} that
$Poly_{d}$ always includes the factorized polynomial of the form,
\begin{equation}
\prod_{j=1}^{d-1}(\frac{(d-j)x+jy}{d}+z_{j}).
\label{fund}
\end{equation}
We try to interpret this formula combinatorially. First, we prepare the
array of the following $d$ circles.
\begin{equation}
\overbrace{\circ\;\circ\;\circ\;\cdots\;\circ\;\circ}^{\mbox{$d$ circles}}
\end{equation}
 Then we can see $d-1$ gaps, and we label 
these gaps as $j \;(j=1,\cdots,d-1)$ from the left to the right.
In this set up, we can regard that the formula (\ref{fund}) 
is summing up all the
configurations whether we insert $z_{j}$ or ${\frac{(d-j)x+jy}{d}}$ into
the $j$-th gap. 
A graphical example is given by the following figure,
\begin{equation}
\circ\;|\;\circ\;|\;\circ\;\circ\Longleftrightarrow
z_{1}z_{2}(\frac{x+3y}{4}),
\end{equation}
where we insert $z_{1}$ and $z_{2}$ to the first and the second gaps 
respectively, and put $(\frac{x+3y}{4})$ in the third gap.
With this interpretation, we search for the extension of 
(\ref{fund}) that includes all the terms appearing in $Poly_{d}$. 
Generally, we cannot expect the factorization of $Poly_{d}$ like 
(\ref{fund}). Instead, we introduce the $d-1$ $U(1)$ charge variables 
$t_{1},t_{2},\cdots,t_{d-1}$ and expect that factorization 
characteristics is conserved by including the $d-1$ $U(1)$ charges. But
afterwards, we have to pick up the neutral charge portion from 
the factorized formula, which is no longer factorized. 
To see if these expectations hold, we carefully look at the $d=4$
case. The difference between $Poly_{4}$ and (\ref{fund}) is given by,   
\begin{eqnarray}
&&Poly_{4}-
(\frac{3x+y}{4}+z_{1})(\frac{x+y}{2}+z_{2})(\frac{x+3y}{4}+z_{3})\no\\
&&=(\frac{3x+11y}{18})(z_{1})^{2}
           +\frac{2}{9}(z_{1})^{3}+(\frac{3x+3y}{8})(z_{2})^{2} 
     +\frac{1}{4}(z_{2})^{3}+(\frac{11x+3y}{18})(z_{3})^{2}
    +\frac{2}{9}(z_{3})^{3}\nonumber\\
  &&+\frac{1}{2}(z_{2})^{2}z_{3}+\frac{1}{3}z_{2}(z_{3})^{2}
  +\frac{2}{3}(z_{1})^{2}z_{3}+\frac{2}{3}z_{1}(z_{3})^{2}
+\frac{1}{3}(z_{1})^{2}z_{2}+\frac{1}{2}z_{1}(z_{2})^{2}.
\label{chigai}
\end{eqnarray}
This formula tells us that we have to admit the insertion of 
$(z_{j})^{n} \;(n=1,2,3)$ into the $j$-th gap. Moreover, we 
found the following equalities:
\begin{eqnarray}
\frac{3x+11y}{18}&=&
\frac{2}{3}\cdot(\frac{x+3y}{4})+\frac{1}{3}\cdot(\frac{x+y}{2})-
\frac{2}{3}\cdot\frac{1}{3}\cdot\frac{3x+y}{4}\no\\
\frac{3x+3y}{8}&=&
\frac{1}{2}\cdot(\frac{x+3y}{4})+\frac{1}{2}\cdot(\frac{3x+y}{4})-
\frac{1}{2}\cdot\frac{1}{2}\cdot\frac{x+y}{2}\no\\
\frac{2}{9}&=&\frac{1}{3}\cdot\frac{2}{3}\qquad
\frac{1}{4}=\frac{1}{2}\cdot\frac{1}{2}.
\label{find}
\end{eqnarray}
Thus, we speculate that the minimal factorization polynomial 
that includes $Poly_{4}$ as the neutral charge part is given by,
\begin{eqnarray}
&&(\frac{3x+y}{4}+z_1+\frac{(z_{1})^{2}}{t_1}(1-\frac{3x+y}{4t_{1}})+
\frac{(z_{1})^{3}}{(t_1)^{2}}+\frac{1}{2}t_{2}+\frac{1}{3}t_{3})\times
\nonumber\\
&&(\frac{2x+2y}{4}+z_2+\frac{(z_{2})^{2}}{t_2}(1-\frac{2x+2y}{4t_{2}})+
   \frac{(z_{2})^{3}}{(t_2)^{2}}+\frac{2}{3}t_{1}+\frac{2}{3}t_{3})\times
\nonumber\\
&&(\frac{x+3y}{4}+z_3+\frac{(z_{3})^{2}}{t_3}(1-\frac{x+3y}{4t_{3}})+
   \frac{(z_{3})^{3}}{(t_3)^{2}}+\frac{1}{3}t_{1}+\frac{1}{2}t_{2}).
\label{yon}
\end{eqnarray} 
Using the same observation, we first expect that the desired factorized
polynomial including $t_{1},\cdots,t_{d-1}$ in the $d=5$ case is given
by,
\begin{eqnarray}
&&\biggl(\frac{4x+y}{5}+\frac{1}{2}t_{2}+\frac{1}{3}t_{3}
+\frac{1}{4}t_{4}+z_1+\no\\
&&\frac{(z_{1})^{2}}{t_1}
(1-(\frac{4x+y}{5t_{1}})+(\frac{4x+y}{5t_{1}})^{2})+
\frac{(z_{1})^{3}}{(t_1)^{2}}
(1-2(\frac{4x+y}{5t_{1}}))+\frac{(z_{1})^{4}}{(t_1)^{3}}\biggr)
\times\nonumber\\
&&\biggl(\frac{3x+2y}{5}+\frac{3}{4}t_{1}+\frac{2}{3}t_{3}+\frac{1}{2}t_{4}+
z_2+\nonumber\\
&&\frac{(z_{2})^{2}}{t_2}
(1-(\frac{3x+2y}{5t_{2}})+(\frac{3x+2y}{5t_{2}})^{2}) +
\frac{(z_{2})^{3}}{(t_2)^{2}}
(1-2(\frac{3x+2y}{5t_{2}}))+\frac{(z_{2})^{4}}{(t_2)^{3}}\biggr)
\times\nonumber\\
&&\biggl(\frac{2x+3y}{5}+\frac{1}{2}t_{1}+\frac{2}{3}t_{2}+\frac{3}{4}t_{4}+
z_3+\nonumber\\
&&\frac{(z_{3})^{2}}{t_3}(1-(\frac{2x+3y}{5t_{3}})+(\frac{2x+3y}{5t_{3}})^{2})+
\frac{(z_{3})^{3}}{(t_3)^{2}}
(1-2(\frac{2x+3y}{5t_{3}}))+\frac{(z_{3})^{4}}{(t_3)^{3}}\biggr)
\times\nonumber\\
&&\biggl(\frac{x+4y}{5}+\frac{1}{4}t_{1}+\frac{1}{3}t_{2}+\frac{1}{2}t_{3}+
z_4+\nonumber\\
&&\frac{(z_{4})^{2}}{t_4}(1-(\frac{x+4y}{5t_{4}})+(\frac{x+4y}{5t_{4}})^{2})+
\frac{(z_{4})^{3}}{(t_4)^{2}}
(1-2(\frac{x+4y}{5t_{4}}))+\frac{(z_{4})^{4}}{(t_4)^{3}}\biggr).
\label{go}
\end{eqnarray}
Looking at (\ref{yon}) and (\ref{go}), we can imagine that the terms 
including $z_{j}$ come from picking up the  
relevant terms of the universal rational function,
\begin{eqnarray} 
&&\sum_{n=0}^{\infty}\frac{(z_{j})^{n+1}}{(t_{j})^{n}}
\sum_{m=0}^{\infty}{-n\choose m}(\frac{(d-j)x+jy}{d t_{j}})^{m}\no\\
&=&z_{j}(t_{j}+\frac{(d-j)x+jy}{d})/
(t_{j}+\frac{(d-j)x+jy}{d}-z_{j}).
\end{eqnarray}
Hence our first trial of the factorized formula including $Poly_{d}$
as the neutral charge portion takes the following form,
\begin{eqnarray}
&&\prod_{j=1}^{d-1}
\biggl(\frac{(d-j)x+jy}{d}+\sum_{i=1}^{j-1}\frac{d-j}{d-i}t_{i}+
\sum_{i=j+1}^{d-1}\frac{j}{i}t_{i}+\no\\
&&z_{j}(t_{j}+\frac{(d-j)x+jy}{d})/(t_{j}+\frac{(d-j)x+jy}{d}-z_{j})
\biggr).
\label{trial1} 
\end{eqnarray}
In the $d=5$ case, (\ref{trial1}) reproduces almost all the terms in the 
$d=5$ recursive formula determined in \cite{cj}
except for the coefficients of $(z_{i})^{2}(z_{j})^{2}\;(i<j)$.

Therefore, let us  consider the exceptional terms.
For example, we pick up the coefficient of the monomial 
$(z_{1})^{2}(z_{4})^{2}$, which corresponds to the term : 
\begin{equation}
L^{N+1,k,1}_{n}L^{N+1,k,3}_{n-1+N-k}L^{N+1,k,1}_{n+4(N-k)}.
\end{equation} 
The prediction of the coefficient
of $(z_{1})^{2}(z_{4})^{2}$ by (\ref{trial1}) is 
\begin{equation}
\frac{3}{4}\cdot\frac{3}{4}+\frac{1}{2}\cdot\frac{1}{2}=\frac{13}{16}.
\end{equation}
But the true coefficient of $(z_{1})^{2}(z_{4})^{2}$ is $\frac{2}{3}$
according to Appendix A. So we have to modify the formula (\ref{trial1}) to 
reproduce the correct coefficient $\frac{2}{3}$. 

At this stage, we look back at the formula (\ref{trial1}) and pay attention to the 
term:
\begin{equation}
\frac{(d-j)x+jy}{d}+\sum_{i=1}^{j-1}\frac{d-j}{d-i}t_{i}+
\sum_{i=j+1}^{d-1}\frac{j}{i}t_{i}.
\label{att}
\end{equation}
 In (\ref{att}), we can interpret formally $x$ and $y$ as $t_{0}$
and $t_{d}$ respectively. Moreover, it is natural to add $t_{j}$ 
to (\ref{att}), which is irrelevant in picking up the neutral charge
portion. In sum, we guess that if the term $\frac{(d-j)x+jy}{d}$ appear 
in the factorized formula, it is always accompanied by 
$t_{j}+\sum_{i=1}^{j-1}\frac{d-j}{d-i}t_{i}+
\sum_{i=j+1}^{d-1}\frac{j}{i}t_{i}$. In this way,
 we modify (\ref{trial1}) into 
the form:
\begin{eqnarray}
&&\prod_{j=1}^{d-1}
\biggl(\frac{(d-j)x+jy}{d}+\sum_{i=1}^{j}\frac{d-j}{d-i}t_{i}+
\sum_{i=j+1}^{d-1}\frac{j}{i}t_{i}+\no\\
&&z_{j}(\frac{(d-j)x+jy}{d}+\sum_{i=1}^{j}\frac{d-j}{d-i}t_{i}+
\sum_{i=j+1}^{d-1}\frac{j}{i}t_{i})/
(\frac{(d-j)x+jy}{d}+\sum_{i=1}^{j}\frac{d-j}{d-i}t_{i}+
\sum_{i=j+1}^{d-1}\frac{j}{i}t_{i}-z_{j})\biggr),\no\\
\label{trial2} 
\end{eqnarray}
that is the formula given in (\ref{p}).
With the formula (\ref{trial2}), we compute the coefficient of
$(z_{1})^{2}(z_{4})^{2}$. If we expand the $j$-th factor of (\ref{trial2}) 
in powers of $\frac{1}{t_{j}}$, we can see, in the case of
$(z_{1})^{2}(z_{4})^{2}$, that we have infinite contributions 
from the terms like $\frac{(t_{1})^{n}}{(t_{4})^{n+1}}\cdot  t_{1}\cdot t_{4}
\cdot\frac{(t_{4})^{n}}{(t_{1})^{n+1}}\;(n\geq 0)$, whose $j$-th factor
$(j=1,2,3,4)$ comes from the $j$-th factor of (\ref{trial2}) respectively.
And we have to sum up these contributions. This operation 
is done as follows:
\begin{eqnarray}
&&(\mbox{coefficients of $(z_{1})^{2}(z_{4})^{2}$})\no\\
&&=\biggl(\mbox{constant term of}\;(\frac{3}{4}t_{1}+\frac{1}{2}t_{4})
\cdot(\frac{1}{2}t_{1}+\frac{3}{4}t_{4})\cdot
(\frac{1}{t_{1}}\sum_{i=0}^{\infty}(-1)^{i}(\frac{t_{4}}{4t_{1}})^{i})\cdot
(\frac{1}{t_{4}}\sum_{j=0}^{\infty}(-1)^{j}(\frac{t_{1}}{4t_{4}})^{j})
\biggr)\no\\
&&=(\frac{3}{4}\cdot\frac{3}{4}+\frac{1}{2}\cdot\frac{1}{2})\cdot\frac{16}{15}
-(\frac{1}{2}\cdot\frac{3}{4}+\frac{3}{4}\cdot\frac{1}{2})\cdot\frac{4}{15}
=\frac{2}{3}.
\end{eqnarray}
Thus, we can see that (\ref{trial2}) reproduces the right coefficient of 
$(z_{1})^{2}(z_{4})^{2}$! Similar computations lead us to complete
agreement between the prediction of (\ref{trial2}) and the $d=5$
recursive formula in Appendix A.    
\section{Test for the $d=6$ case}
Now, we test our main conjecture for the $d=6$ rational curves.
First, we write down the relevant part of $R_{6}$:
\begin{eqnarray}
&&R_{6}=\no\\
&&\biggl(\frac{5x+y}{6}+\frac{1}{2}t_{2}+\frac{1}{3}t_{3}
+\frac{1}{4}t_{4}+\frac{1}{5}t_{5}+z_{1}+
\frac{(z_{1})^{2}}{t_1}(t_{1}/(t_{1}+\frac{5x+y}{6}+
\frac{1}{2}t_{2}+\frac{1}{3}t_{3}+\frac{1}{4}t_{4}+\frac{1}{5}t_{5}))+\no\\
&&\frac{(z_{1})^{3}}{(t_1)^{2}}(t_{1}/(t_{1}+\frac{5x+y}{6}+
\frac{1}{2}t_{2}+\frac{1}{3}t_{3}+\frac{1}{4}t_{4}+\frac{1}{5}t_{5}))^{2}
+\frac{(z_{1})^{4}}{(t_1)^{3}}(1-3(\frac{5x+y}{6t_{1}}))+
\frac{(z_{1})^{5}}{(t_1)^{4}})\biggr)
\times\no\\
&&\biggl(\frac{4x+2y}{6}+\frac{4}{5}t_{1}+\frac{2}{3}t_{3}
+\frac{2}{4}t_{4}+\frac{2}{5}t_{5}+z_{2}+
\frac{(z_{2})^{2}}{t_{2}}(t_{2}/(t_{2}+\frac{4x+2y}{6}+
\frac{4}{5}t_{1}+\frac{2}{3}t_{3}+\frac{2}{4}t_{4}+\frac{2}{5}t_{5}))+\no\\
&&\frac{(z_{2})^{3}}{(t_{2})^{2}}(t_{2}/(t_{2}+\frac{4x+2y}{6}+
\frac{4}{5}t_{1}+\frac{2}{3}t_{3}+\frac{2}{4}t_{4}+\frac{2}{5}t_{5}))^{2}
+\frac{(z_{2})^{4}}{(t_{2})^{3}}(1-3(\frac{4x+2y}{6t_{2}}))+
\frac{(z_{2})^{5}}{(t_{2})^{4}})\biggr)
\times\no\\
&&\biggl(\frac{3x+3y}{6}+\frac{3}{5}t_{1}+\frac{3}{4}t_{2}
+\frac{3}{4}t_{4}+\frac{3}{5}t_{5}+z_{3}+
\frac{(z_{3})^{2}}{t_{3}}(t_{3}/(t_{3}+\frac{3x+3y}{6}+
\frac{3}{5}t_{1}+\frac{3}{4}t_{2}+\frac{3}{4}t_{4}+\frac{3}{5}t_{5}))+\no\\
&&\frac{(z_{3})^{3}}{(t_{3})^{2}}(t_{3}/(t_{3}+\frac{3x+3y}{6}+
\frac{3}{5}t_{1}+\frac{3}{4}t_{2}+\frac{3}{4}t_{4}+\frac{3}{5}t_{5}))^{2}
+\frac{(z_{3})^{4}}{(t_{3})^{3}}(1-3(\frac{3x+3y}{6t_{3}}))+
\frac{(z_{3})^{5}}{(t_{3})^{4}})\biggr)
\times\no\\
&&\biggl(\frac{2x+4y}{6}+\frac{2}{5}t_{1}+\frac{2}{4}t_{2}
+\frac{2}{3}t_{3}+\frac{4}{5}t_{5}+z_{4}+
\frac{(z_{4})^{2}}{t_{4}}(t_{4}/(t_{4}+\frac{2x+4y}{6}+
\frac{2}{5}t_{1}+\frac{2}{4}t_{2}+\frac{2}{3}t_{3}+\frac{4}{5}t_{5}))+\no\\
&&\frac{(z_{4})^{3}}{(t_{4})^{2}}(t_{4}/(t_{4}+\frac{2x+4y}{6}+
\frac{2}{5}t_{1}+\frac{2}{4}t_{2}+\frac{2}{3}t_{3}+\frac{4}{5}t_{5}))^{2}
+\frac{(z_{4})^{4}}{(t_{4})^{3}}(1-3(\frac{2x+4y}{6t_{4}}))+
\frac{(z_{4})^{5}}{(t_{4})^{4}})\biggr)
\times\no\\
&&\biggl(\frac{x+5y}{6}+\frac{1}{5}t_{1}+\frac{1}{4}t_{2}
+\frac{1}{3}t_{3}+\frac{1}{2}t_{4}+z_{5}+
\frac{(z_{5})^{2}}{t_{5}}(t_{5}/(t_{5}+\frac{x+5y}{6}+
\frac{1}{5}t_{1}+\frac{1}{4}t_{2}+\frac{1}{3}t_{3}+\frac{1}{2}t_{4}))+\no\\
&&\frac{(z_{5})^{3}}{(t_{5})^{2}}(t_{5}/(t_{5}+\frac{x+5y}{6}+
\frac{1}{5}t_{1}+\frac{1}{4}t_{2}+\frac{1}{3}t_{3}+\frac{1}{2}t_{5}))^{2}
+\frac{(z_{5})^{4}}{(t_{5})^{3}}(1-3(\frac{x+5y}{6t_{5}}))+
\frac{(z_{5})^{5}}{(t_{5})^{4}})\biggr)\no\\
\end{eqnarray}
Like the $d=5$ cases, we have to evaluate some non-trivial residue 
integrals of some rational functions in $t_{1},\cdots,t_{5}$.
One of the typical examples appears in evaluating the coefficient of 
$(z_{1})^{3}(z_{2})^{2}$:
\begin{eqnarray}
\frac{1}{(2\pi\sqrt{-1})^{2}}\int_{C_{1}}dt_{1}
\int_{C_{2}}dt_{2}
\frac{(\frac{3}{5}t_{1}+\frac{3}{4}t_{2})
(\frac{2}{5}t_{1}+\frac{2}{4}t_{2})
(\frac{1}{5}t_{1}+\frac{1}{4}t_{2})}
{t_{1}t_{2}(t_{1}+
\frac{1}{2}t_{2})^{2}(t_{2}+
\frac{4}{5}t_{1})}=\frac{3}{50}.
\label{period}
\end{eqnarray}
 After the tedious but elementary 
calculation like (\ref{period}), we obtain the following formula of $Poly_{6}$.
\begin{eqnarray}
Poly_{6}&=&(5/324)x^{5}+(29/216)x^{4}y+(227/648)x^{3}y^{2}+
(227/648)x^{2}y^{3}+(29/216)xy^{4}+(5/324)y^{5}\nonumber\\
&&\nonumber\\
&&+z_{1}((1/54)x^{4}+(17/108)x^{3}y+(7/18)x^{2}y^{2}+(37/108)xy^{3}+(5/54)y^{4}
\nonumber\\
&&+((1/45)x^{3}+(83/450)x^{2}y+(967/2250)xy^{2}+(1829/5625)y^{3})z_{1}
\nonumber\\
&&+((2/75)x^{2}+(27/125)xy+(886/1875)y^{2})(z_{1})^{2}\nonumber\\
&&+((4/125)x+(158/625)y)(z_{1})^{3}+(24/625)(z_{1})^{4})\nonumber\\
&&\nonumber\\
&&+z_{2}((5/216)x^{4}+(41/216)x^{3}y+(31/72)x^{2}y^{2}+
(67/216)xy^{3}+(5/108)y^{4}
\nonumber\\
&&+((5/144)x^{3}+(77/288)x^{2}y+(295/576)xy^{2}+(5/36)y^{3})z_{2}
\nonumber\\
&&+((5/96)x^{2}+(3/8)xy+(53/192)y^{2})(z_{2})^{2}\nonumber\\
&&+((5/64)x+(7/32)y)(z_{2})^{3}+(3/64)(z_{2})^{4})\nonumber\\
&&\nonumber\\
&&+z_{3}((5/162)x^{4}+(77/324)x^{3}y+(25/54)x^{2}y^{2}+
(77/324)xy^{3}+(5/162)y^{4}
\nonumber\\
&&+((5/81)x^{3}+(67/162)x^{2}y+(67/162)xy^{2}+(5/81)y^{3})z_{3}
\nonumber\\
&&+((10/81)x^{2}+(37/81)xy+(10/81)y^{2})(z_{3})^{2}\nonumber\\
&&+((4/27)x+(4/27)y)(z_{3})^{3}+(4/81)(z_{3})^{4})\nonumber\\
&&\nonumber\\
&&+z_{4}((5/108)x^{4}+(67/216)x^{3}y+(31/72)x^{2}y^{2}+
(41/216)xy^{3}+(5/216)y^{4}
\nonumber\\
&&+((5/36)x^{3}+(295/576)x^{2}y+(77/288)xy^{2}+(5/144)y^{3})z_{4}
\nonumber\\
&&+((53/192)x^{2}+(3/8)xy+(5/96)y^{2})(z_{4})^{2}\nonumber\\
&&+((7/32)x+(5/64)y)(z_{4})^{3}+(3/64)(z_{4})^{4})\nonumber\\
&&\nonumber\\
&&+z_{5}((5/54)x^{4}+(37/108)x^{3}y+(7/18)x^{2}y^{2}+
(17/108)xy^{3}+(1/54)y^{4}
\nonumber\\
&&+((1829/5625)x^{3}+(967/2250)x^{2}y+(83/450)xy^{2}+(1/45)y^{3})z_{5}
\nonumber\\
&&+((886/1875)x^{2}+(27/125)xy+(2/75)y^{2})(z_{5})^{2}\nonumber\\
&&+((158/625)x+(4/125)y)(z_{5})^{3}+(24/625)(z_{5})^{4})\nonumber\\
&&\nonumber\\
&&+z_{1}z_{2}
(((1/36)x^{3}+(2/9)x^{2}y+(17/36)xy^{2}+(5/18)y^{3})
\nonumber\\
&&+((1/30)x^{2}+(13/50)xy+(193/375)y^{2})z_{1}
+((1/24)x^{2}+(5/16)xy+(53/96)y^{2})z_{2}
\nonumber\\
&&+((1/25)x+(38/125)y)(z_{1})^{2}
+((1/16)x+(7/16)y)(z_{2})^{2}
+((1/20)x+(73/200)y)z_{1}z_{2}\nonumber\\
&&+(3/50)(z_{1})^{2}z_{2}+(3/40)z_{1}(z_{2})^{2}+
(6/125)(z_{1})^{3}+(3/32)(z_{2})^{3} )\nonumber\\
&&\nonumber\\
&&+z_{1}z_{5}
(((1/9)x^{3}+(7/18)x^{2}y+(7/18)xy^{2}+(1/9)y^{3})
\nonumber\\
&&+((2/15)x^{2}+(11/25)xy+(142/375)y^{2})z_{1}
+((142/375)x^{2}+(11/25)xy+(2/15)y^{2})z_{5}
\nonumber\\
&&+((4/25)x+(62/125)y)(z_{1})^{2}
+((62/125)x+(4/25)y)(z_{5})^{2}
+((11/25)x+(11/25)y)z_{1}z_{5}\nonumber\\
&&+(24/125)(z_{1})^{3}+(24/125)(z_{5})^{3}
+(51/100)(z_{1})^{2}z_{5}+(51/100)z_{1}(z_{5})^{2})\nonumber\\
&&\nonumber\\
&&+z_{4}z_{5}
(((5/18)x^{3}+(17/36)x^{2}y+(2/9)xy^{2}+(1/36)y^{3})
\nonumber\\
&&+((53/96)x^{2}+(5/16)xy+(1/24)y^{2})z_{4}
+((193/375)x^{2}+(13/50)xy+(1/30)y^{2})z_{5}
\nonumber\\
&&+((7/16)x+(1/16)y)(z_{4})^{2}
+((38/125)x+(1/25)y)(z_{5})^{2}
+((73/200)x+(1/20)y)z_{4}z_{5}\nonumber\\
&&+(3/32)(z_{4})^{3}+(6/125)(z_{5})^{3}
+(3/40)(z_{4})^{2}z_{5}+(3/50)z_{4}(z_{5})^{2})\nonumber\\
&&\nonumber\\
&&+z_{1}z_{3}
(((1/27)x^{3}+(5/18)x^{2}y+(1/2)xy^{2}+(5/27)y^{3})
\nonumber\\
&&+((2/45)x^{2}+(73/225)xy+(602/1125)y^{2})z_{1}
+((2/27)x^{2}+(13/27)xy+(10/27)y^{2})z_{3}
\nonumber\\
&&+((4/75)x+(142/375)y)(z_{1})^{2}
+((4/27)x+(4/9)y)(z_{3})^{2}
+((4/45)x+(14/25)y)z_{1}z_{3}\nonumber\\
&&+(8/125)(z_{1})^{3}+(4/27)(z_{3})^{3}
+(8/75)(z_{1})^{2}z_{3}+(8/45)z_{1}(z_{3})^{2})\nonumber\\
&&\nonumber\\
&&+z_{1}z_{4}
(((1/18)x^{3}+(13/36)x^{2}y+(4/9)xy^{2}+(5/36)y^{3})
\nonumber\\
&&+((1/15)x^{2}+(21/50)xy+(337/750)y^{2})z_{1}
+((1/6)x^{2}+(9/16)xy+(5/24)y^{2})z_{4}
\nonumber\\
&&+((2/25)x+(61/125)y)(z_{1})^{2}
+((5/16)x+(5/16)y)(z_{4})^{2}
+((1/5)x+(61/100)y)z_{1}z_{4}\nonumber\\
&&+(12/125)(z_{1})^{3}+(3/16)(z_{4})^{3}
+(6/25)(z_{1})^{2}z_{4}+(7/20)z_{1}(z_{4})^{2})\nonumber\\
&&\nonumber\\
&&+z_{3}z_{4}
(((5/54)x^{3}+(19/36)x^{2}y+(1/3)xy^{2}+(5/108)y^{3})
\nonumber\\
&&+((5/27)x^{2}+(31/54)xy+(5/54)y^{2})z_{3}
+((5/18)x^{2}+(67/144)xy+(5/72)y^{2})z_{4}
\nonumber\\
&&+((2/9)x+(5/27)y)(z_{3})^{2}
+((13/48)x+(5/48)y)(z_{4})^{2}
+((1/3)x+(5/36)y)z_{3}z_{4}\nonumber\\
&&+(2/27)(z_{3})^{3}+(1/16)(z_{4})^{3}
+(1/9)(z_{3})^{2}z_{4}+(1/12)z_{3}(z_{4})^{2})\nonumber\\
&&\nonumber\\
&&+z_{3}z_{5}
(((5/27)x^{3}+(1/2)x^{2}y+(5/18)xy^{2}+(1/27)y^{3})
\nonumber\\
&&+((10/27)x^{2}+(13/27)xy+(2/27)y^{2})z_{3}
+((602/1125)x^{2}+(73/225)xy+(2/45)y^{2})z_{5}
\nonumber\\
&&+((4/9)x+(4/27)y)(z_{3})^{2}
+((142/375)x+(4/75)y)(z_{5})^{2}
+((14/25)x+(4/45)y)z_{3}z_{5}\nonumber\\
&&+(4/27)(z_{3})^{3}+(8/125)(z_{5})^{3}
+(8/45)(z_{3})^{2}z_{5}+(8/75)z_{3}(z_{5})^{2})\nonumber\\
&&\nonumber\\
&&+z_{2}z_{5}
(((5/36)x^{3}+(4/9)x^{2}y+(13/36)xy^{2}+(1/18)y^{3})
\nonumber\\
&&+((5/24)x^{2}+(9/16)xy+(1/6)y^{2})z_{2}
+((337/750)x^{2}+(21/50)xy+(1/15)y^{2})z_{5}
\nonumber\\
&&+((5/16)x+(5/16)y)(z_{2})^{2}
+((61/125)x+(2/25)y)(z_{5})^{2}
+((61/100)x+(1/5)y)z_{2}z_{5}\nonumber\\
&&+(3/16)(z_{2})^{3}+(12/125)(z_{5})^{3}
+(7/20)(z_{2})^{2}z_{5}+(6/25)z_{2}(z_{5})^{2})\nonumber\\
&&\nonumber\\
&&+z_{2}z_{3}
(((5/108)x^{3}+(1/3)x^{2}y+(19/36)xy^{2}+(5/54)y^{3})
\nonumber\\
&&+((5/72)x^{2}+(67/144)xy+(5/18)y^{2})z_{2}
+((5/54)x^{2}+(31/54)xy+(5/27)y^{2})z_{3}
\nonumber\\
&&+((5/48)x+(13/48)y)(z_{2})^{2}
+((5/27)x+(2/9)y)(z_{3})^{2}
+((5/36)x+(1/3)y)z_{2}z_{3}\nonumber\\
&&+(1/16)(z_{2})^{3}+(2/27)(z_{3})^{3}
+(1/12)(z_{2})^{2}z_{3}+(1/9)z_{2}(z_{3})^{2})\nonumber\\
&&\nonumber\\
&&+z_{2}z_{4}
(((5/72)x^{3}+(31/72)x^{2}y+(31/72)xy^{2}+(5/72)y^{3})
\nonumber\\
&&+((5/48)x^{2}+(19/32)xy+(5/24)y^{2})z_{2}
+((5/24)x^{2}+(19/32)xy+(5/48)y^{2})z_{4}
\nonumber\\
&&+((5/32)x+(11/32)y)(z_{2})^{2}
+((11/32)x+(5/32)y)(z_{4})^{2}
+((5/16)x+(5/16)y)z_{2}z_{4}\nonumber\\
&&+(3/32)(z_{2})^{3}+(3/32)(z_{4})^{3}
+(3/16)(z_{2})^{2}z_{4}+(3/16)z_{2}(z_{4})^{2})\nonumber\\
&&\nonumber\\
&&+z_{1}z_{2}z_{3}
(((1/18)x^{2}+(7/18)xy+(5/9)y^{2})+((1/15)x+(34/75)y)z_{1}\no\\
&&+((1/12)x+(13/24)y)z_{2}+((1/9)x+(2/3)y)z_{3}\nonumber\\
&&+(2/25)(z_{1})^{2}+(1/8)(z_{2})^{2}+
(2/9)(z_{3})^{2}+(1/10)z_{1}z_{2}+(1/6)z_{2}z_{3}
+(2/15)z_{1}z_{3})\nonumber\\
&&\nonumber\\
&&+z_{1}z_{2}z_{5}
(((1/6)x^{2}+(1/2)xy+(1/3)y^{2})+((1/5)x+(14/25)y)z_{1}\no\\
&&+((1/4)x+(5/8)y)z_{2}+((13/25)x+(2/5)y)z_{5}\nonumber\\
&&+(6/25)(z_{1})^{2}+(3/8)(z_{2})^{2}+
(12/25)(z_{5})^{2}+(3/10)z_{1}z_{2}+(7/10)z_{2}z_{5}
+(3/5)z_{1}z_{5})\nonumber\\
&&\nonumber\\
&&+z_{1}z_{4}z_{5}
(((1/3)x^{2}+(1/2)xy+(1/6)y^{2})+((2/5)x+(13/25)y)z_{1}\no\\
&&+((5/8)x+(1/4)y)z_{4}+((14/25)x+(1/5)y)z_{5}\nonumber\\
&&+(12/25)(z_{1})^{2}+(3/8)(z_{4})^{2}+
(6/25)(z_{5})^{2}+(7/10)z_{1}z_{4}+(3/10)z_{4}z_{5}
+(3/5)z_{1}z_{5})\nonumber\\
&&\nonumber\\
&&+z_{3}z_{4}z_{5}
(((5/9)x^{2}+(7/18)xy+(1/18)y^{2})+((2/3)x+(1/9)y)z_{3}\no\\
&&+((13/24)x+(1/12)y)z_{4}+((34/75)x+(1/15)y)z_{5}\nonumber\\
&&+(2/9)(z_{3})^{2}+(1/8)(z_{4})^{2}+
(2/25)(z_{5})^{2}+(1/6)z_{3}z_{4}+(1/10)z_{4}z_{5}
+(2/15)z_{3}z_{5})\nonumber\\
&&\nonumber\\
&&+z_{1}z_{2}z_{4}
(((1/12)x^{2}+(1/2)xy+(5/12)y^{2})+((1/10)x+(29/50)y)z_{1}\no\\
&&+((1/8)x+(11/16)y)z_{2}+((1/4)x+(5/8)y)z_{4}\nonumber\\
&&+(3/25)(z_{1})^{2}+(3/16)(z_{2})^{2}+
(3/8)(z_{4})^{2}+(3/20)z_{1}z_{2}+(3/8)z_{2}z_{4}
+(3/10)z_{1}z_{4})\nonumber\\
&&\nonumber\\
&&+z_{1}z_{3}z_{4}
(((1/9)x^{2}+(11/18)xy+(5/18)y^{2})+((2/15)x+(53/75)y)z_{1}\no\\
&&+((2/9)x+(5/9)y)z_{3}+((1/3)x+(5/12)y)z_{4}\nonumber\\
&&+(4/25)(z_{1})^{2}+(2/9)(z_{3})^{2}+
(1/4)(z_{4})^{2}+(4/15)z_{1}z_{3}+(1/3)z_{3}z_{4}
+(2/5)z_{1}z_{4})\nonumber\\
&&\nonumber\\
&&+z_{2}z_{3}z_{4}
(((5/36)x^{2}+(13/18)xy+(5/36)y^{2})+((5/24)x+(5/12)y)z_{2}\no\\
&&+((5/18)x+(5/18)y)z_{3}+((5/12)x+(5/24)y)z_{4}\nonumber\\
&&+(1/8)(z_{2})^{2}+(1/9)(z_{3})^{2}+
(1/8)(z_{4})^{2}+(1/6)z_{2}z_{3}+(1/4)z_{2}z_{4}
+(1/6)z_{3}z_{4})\nonumber\\
&&\nonumber\\
&&+z_{1}z_{3}z_{5}
(((2/9)x^{2}+(5/9)xy+(2/9)y^{2})+((4/15)x+(46/75)y)z_{1}\no\\
&&+((4/9)x+(4/9)y)z_{3}+((46/75)x+(4/15)y)z_{5}\nonumber\\
&&+(8/25)(z_{1})^{2}+(4/9)(z_{3})^{2}+
(8/25)(z_{5})^{2}+(8/15)z_{1}z_{3}+(8/15)z_{3}z_{5}
+(7/10)z_{1}z_{5})\nonumber\\
&&\nonumber\\
&&+z_{2}z_{3}z_{5}
(((5/18)x^{2}+(11/18)xy+(1/9)y^{2})+((5/12)x+(1/3)y)z_{2}\no\\
&&+((5/9)x+(2/9)y)z_{3}+((53/75)x+(2/15)y)z_{5}\nonumber\\
&&+(1/4)(z_{2})^{2}+(2/9)(z_{3})^{2}+
(4/25)(z_{5})^{2}+(1/3)z_{2}z_{3}+(4/15)z_{3}z_{5}
+(2/5)z_{2}z_{5})\nonumber\\
&&\nonumber\\
&&+z_{2}z_{4}z_{5}
(((5/12)x^{2}+(1/2)xy+(1/12)y^{2})+((5/8)x+(1/4)y)z_{2}\no\\
&&+((11/16)x+(1/8)y)z_{4}+((29/50)x+(1/10)y)z_{5}\nonumber\\
&&+(3/8)(z_{2})^{2}+(3/16)(z_{4})^{2}+
(3/25)(z_{5})^{2}+(3/8)z_{2}z_{4}+(3/20)z_{4}z_{5}
+(3/10)z_{2}z_{5})\nonumber\\
&&\nonumber\\
&&+z_{1}z_{2}z_{3}z_{4}
(((1/6)x+(5/6)y)+(1/5)z_{1}+(1/4)z_{2}+
(1/3)z_{3}+(1/2)z_{4})\nonumber\\
&&\nonumber\\
&&+z_{1}z_{2}z_{3}z_{5}
(((1/3)x+(2/3)y)+(2/5)z_{1}+(1/2)z_{2}+
(2/3)z_{3}+(4/5)z_{5})\nonumber\\
&&\nonumber\\
&&+z_{1}z_{2}z_{4}z_{5}
(((1/2)x+(1/2)y)+(3/5)z_{1}+(3/4)z_{2}+
(3/4)z_{4}+(3/5)z_{5})\nonumber\\
&&\nonumber\\
&&+z_{1}z_{3}z_{4}z_{5}
(((2/3)x+(1/3)y)+(4/5)z_{1}+(2/3)z_{3}+
(1/2)z_{4}+(2/5)z_{5})\nonumber\\
&&\nonumber\\
&&+z_{2}z_{3}z_{4}z_{5}
(((5/6)x+(1/6)y)+(1/2)z_{2}+(1/3)z_{3}+
(1/4)z_{4}+(1/5)z_{5})\nonumber\\
&&\nonumber\\
&&+z_{1}z_{2}z_{3}z_{4}z_{5}
\label{6}
\end{eqnarray}
The first non-trivial test for this formula is the compatibility with 
the conjecture proposed in \cite{cj}.
We restate here the conjecture:
\begin{conj}
Formal iteration of the recursive formula of Fano hypersurfaces
 for descending $N$
down to the case $N=k$  yields the coefficients of the hypergeometric
 series used in the mirror calculation, i.e., we should have the 
equalities:
\begin{eqnarray}
a_{d}&=&\tilde{L}_{0}^{k,k,d}\nonumber\\
b_{d}&=&\frac{1}{d}\tilde{L}_{1}^{k,k,d}+
\sum_{m=1}^{d-1}\frac{1}{m}\tilde{L}_{1}^{k,k,m}\cdot\tilde{L}_{0}^{k,k,d-m}.
\label{hyper}
\end{eqnarray}
 Here
$$a_d  = \frac{(kd)!}{(d!)^{k}} \; ,\;\;\;\;
 b_d = a_d(\sum_{i=1}^{d}\sum_{m=1}^{k-1}\frac{m}{i(ki-m)})$$
are the coefficients of the hypergeometric
series associated to the solutions
$$w_{0}(x) = \sum_{d=0}^{\infty}a_{d}e^{dx},\;\;\;\;
w_{1}(x) = \sum_{d=1}^{\infty}b_{d}e^{dx}+w_{0}(x)x $$
of the differential equation
\begin{eqnarray}
&&((\frac{d}{dx})^{k-1}-ke^{x}(k\frac{d}{dx}+1)
(k\frac{d}{dx}+2)\cdots(k\frac{d}{dx}+k-1))w_{i}(x)=0.
\nonumber\\
\label{diff}
\end{eqnarray}
\end{conj}
We tested the recursive formula for the $d=6$ rational curves 
obtained from the formula (\ref{6}) by checking that it reproduces the
coefficients $a_{6}$ and $b_{6}$ for the $k\leq 20$ cases. 
Of course, we argue that 
in the $N>k$ region, this recursive formula predicts the rational 
structural constant $L_{n}^{N,k,6}$ 
for the $d=6$ rational curves. 
For example, we predict 
\begin{equation}
L_{0}^{8,7,6}=13799153353276807722049582771200.
\end{equation}
Finally, we discuss the validity of Conjecture 1 for higher degree 
cases. In these cases, the most non-trivial part of the calculation using 
Conjecture 1 comes from the terms like 
$(z_{1})^{2}(z_{2})^{2}\cdots(z_{m})^{2}$. In evaluating the
coefficients of these terms, we have to integrate the rational 
functions in many variables. For example, we consider the coefficient of 
$(z_{1})^{2}(z_{2})^{2}(z_{3})^{2}$ in the $d=7$ case. This coefficient 
is evaluated by the residue integral of the following rational function:
\begin{eqnarray}
&&\frac{1}{(2\pi\sqrt{-1})^{3}}
\int_{C_{1}}\frac{dt_{1}}{t_{1}}\int_{C_{2}}\frac{dt_{2}}{t_{2}} 
\int_{C_{3}}\frac{dt_{3}}{t_{3}}
\frac{(\frac{3}{6}t_{1}+\frac{3}{5}t_{2}+\frac{3}{4}t_{3})
(\frac{2}{6}t_{1}+\frac{2}{5}t_{2}+\frac{2}{4}t_{3}) 
(\frac{1}{6}t_{1}+\frac{1}{5}t_{2}+\frac{1}{4}t_{3})}
{(t_{1}+\frac{1}{2}t_{2}+\frac{1}{3}t_{3})
(\frac{5}{6}t_{1}+t_{2}+\frac{2}{3}t_{3}) 
(\frac{4}{6}t_{1}+\frac{4}{5}t_{2}+t_{3})}=\frac{1}{20}.\no\\
\label{3}
\end{eqnarray}
The rational number $\frac{1}{20}$ has already appeared as the coefficient 
of $x(z_{1})^{2}(z_{2})^{2}$ 
in the $d=6$ case. This coincidence is compatible with 
the relation (\ref{jincolli}), which suggests many relations between 
the coefficients of $Poly_{d}$'s with different $d$'s. Thus, we can 
rely on our predictions using the residue integral of rational functions 
in many variables like (\ref{3}).   
 We believe that our main conjecture is valid 
in the case of the rational curves of arbitrary degree!
\section{Conclusion}
We first discuss applications of our main conjecture. 
The recursive formulas for descending the dimension of the 
hypersurface while conserving its degree are first introduced by 
A.Collino in an attempt to prove the relation (\ref{jincolli}). Now,
we have the general form of the recursive formula and it
is in principle possible to prove (\ref{jincolli}) along the line 
proposed in \cite{cj} (of course, we have to prove Conjecture 1 before 
this). We also have to consider the rigorous connection between the 
recursive formulas and the hypergeometric series. This line of thought 
also contributes to the analysis of the quantum cohomology ring of the
general type hypersurface along the line of \cite{gene}. 

Physically, 
we can try to interpret our approach as the large $N$ expansion 
of the gauged linear sigma model. From this point of view, we can regard the
recursive formulas as the (discretized) differential equations that 
describe the behavior of the model in the complex $N$-plane. 
Since the formula given in our main conjecture is rather simple,
it might be possible to derive them through the analysis of the gauged 
linear sigma model. 

Geometrically, it seems to be natural to regard the appearance of $d-1$ 
$U(1)$ charge variables $t_{1},t_{2},\cdots ,t_{d-1}$ as the consequence 
of application of some kind of fixed point theorem like the discussion
in \cite{bert}. And it is of course interesting to prove the formulas 
from this line of thought.
          
{\bf Acknowledgment}\\
We would like to thank A.Collino, T.Eguchi S.Hosono and A. Matsuo for
discussions.  We are also grateful to G.Yamamoto for manipulation of 
computers. 
\newpage
\section*{ Appendix A:  Recursive Formulas for Fano Hypersurfaces
with $c_{1}(M_{N}^{k})\geq 2$}
\begin{eqnarray}
L^{N,k,1}_{m}&=&L^{N+1,k,1}_{m}:=L^{k}_{m}\\
L^{N,k,2}_{m}&=&\frac{1}{2}(L^{N+1,k,2}_{m-1}+L^{N+1,k,2}_{m}
+2L^{N+1,k,1}_{m}\cdot L^{N+1,k,1}_{m+(N-k)})\\
L^{N,k,3}_{m}&=&\frac{1}{18}(4L^{N+1,k,3}_{m-2}+10L^{N+1,k,3}_{m-1}
+4L^{N+1,k,3}_{m}\nonumber\\
&&+12L^{N+1,k,2}_{m-1}\cdot L^{N+1,k,1}_{m+2(N-k)}
+9L^{N+1,k,2}_{m}\cdot L^{N+1,k,1}_{m+2(N-k)}\nonumber\\
&&+6L^{N+1,k,2}_{m}\cdot L^{N+1,k,1}_{m+1+2(N-k)}\nonumber\\
&&+6 L^{N+1,k,1}_{m-1}\cdot L^{N+1,k,2}_{m-1+(N-k)}
+9 L^{N+1,k,1}_{m}\cdot L^{N+1,k,2}_{m-1+(N-k)}\nonumber\\
&&+12 L^{N+1,k,1}_{m}\cdot L^{N+1,k,2}_{m+(N-k)}\nonumber\\
&&+18L^{N+1,k,1}_{m}\cdot L^{N+1,k,1}_{m+(N-k)}\cdot L^{N+1,k,1}_{m+2(N-k)})
\label{rec0}
\end{eqnarray}
In the following, we omit the subscripts $N+1$ and $k$ of 
$L_{n}^{N+1,k,4}$ in the r.h.s. 
for brevity.
\begin{eqnarray}
L_{n}^{N,k,4}
&=&\frac{1}{32}(3L_{n-3}^{4}+13L_{n-2}^{4}
+13L_{n-1}^{4}+3L_{n}^{4})\nonumber\\
&&+\frac{1}{72}(9L_{n-2}^{1}L_{n-2+N-k}^{3}+12L_{n-1}^{1}L_{n-2+N-k}^{3}
+16L_{n}^{1}L_{n-2+N-k}^{3}\nonumber\\
&&+36L_{n-1}^{1}L_{n-1+N-k}^{3}
+44L_{n}^{1}L_{n-1+N-k}^{3}+27L_{n}^{1}L_{n+N-k}^{3})\nonumber\\
&&+\frac{1}{16}(3L_{n-2}^{2}L_{n-1+2(N-k)}^{2}
+6L_{n-1}^{2}L_{n-1+2(N-k)}^{2}
+4L_{n}^{2}L_{n-1+2(N-k)}^{2}\nonumber\\
&&+10L_{n-1}^{2}L_{n+2(N-k)}^{2}
+6L_{n}^{2}L_{n+2(N-k)}^{2}+3L_{n}^{2}L_{n+1+2(N-k)}^{2})
\nonumber\\
&&+\frac{1}{72}(27L_{n-2}^{3}L_{n+3(N-k)}^{1}
+44L_{n-1}^{3}L_{n+3(N-k)}^{1}
+16L_{n}^{3}L_{n+3(N-k)}^{1}\nonumber\\
&&+36L_{n-1}^{3}L_{n+1+3(N-k)}^{1}
+12L_{n}^{3}L_{n+1+3(N-k)}^{1}+
9L_{n}^{3}L_{n+2+3(N-k)}^{1})\nonumber\\
&&+\frac{1}{12}(3L_{n-1}^{1}L_{n-1+N-k}^{1}L_{n-1+2(N-k)}^{2}
+4L_{n}^{1}L_{n-1+N-k}^{1}L_{n-1+2(N-k)}^{2}\nonumber\\
&&+6L_{n}^{1}L_{n+N-k}^{1}L_{n-1+2(N-k)}^{2}
+9L_{n}^{1}L_{n+N-k}^{1}L_{n+2(N-k)}^{2})\nonumber\\
&&+\frac{1}{6}(3L_{n-1}^{1}L_{n-1+N-k}^{2}L_{n+3(N-k)}^{1}
+4L_{n}^{1}L_{n-1+N-k}^{2}L_{n+3(N-k)}^{1}\nonumber\\
&&+4L_{n}^{1}L_{n+N-k}^{2}L_{n+3(N-k)}^{1}
+3L_{n}^{1}L_{n+N-k}^{2}L_{n+1+3(N-k)}^{1})\nonumber\\
&&+\frac{1}{12}(9L_{n-1}^{2}L_{n+2(N-k)}^{1}L_{n+3(N-k)}^{1}
+6L_{n}^{2}L_{n+2(N-k)}^{1}L_{n+3(N-k)}^{1}\nonumber\\
&&+4L_{n}^{2}L_{n+1+2(N-k)}^{1}L_{n+3(N-k)}^{1}
+3L_{n}^{2}L_{n+1+2(N-k)}^{1}L_{n+1+3(N-k)}^{1})\nonumber\\
&&+L_{n}^{1}L_{n+N-k}^{1}L_{n+2(N-k)}^{1}L_{n+3(N-k)}^{1}.    
\end{eqnarray}
\newpage
\begin{eqnarray}
&&L^{N,k,5}_{n} \nonumber\\
&=&\frac{24}{625}L_{n-4}^{5}+\frac{154}{625}L_{n-3}^{5}
+\frac{269}{625}L_{n-2}^{5}+\frac{154}{625}L_{n-1}^{5}
+\frac{24}{625}L_{n}^{5} \nonumber\\
&&+\frac{6}{125}L_{n-3}^{1}L_{n-3+N-k}^{4}
+\frac{3}{50}L_{n-2}^{1}L_{n-3+N-k}^{4}
+\frac{3}{40}L_{n-1}^{1}L_{n-3+N-k}^{4}\nonumber\\
&&+\frac{3}{32}L_{n}^{1}L_{n-3+N-k}^{4}
+\frac{37}{125}L_{n-2}^{1}L_{n-2+N-k}^{4}
+\frac{71}{200}L_{n-1}^{1}L_{n-2+N-k}^{4}\nonumber\\
&&+\frac{17}{40}L_{n}^{1}L_{n-2+N-k}^{4}
+\frac{58}{125}L_{n-1}^{1}L_{n-1+N-k}^{4}
+\frac{393}{800}L_{n}^{1}L_{n-1+N-k}^{4}\nonumber\\
&&+\frac{24}{125}L_{n}^{1}L_{n+N-k}^{4}\nonumber\\
&&+\frac{8}{125}L_{n-3}^{2}L_{n-2+2(N-k)}^{3}+
\frac{8}{75}L_{n-2}^{2}L_{n-2+2(N-k)}^{3}
+\frac{8}{45}L_{n-1}^{2}L_{n-2+2(N-k)}^{3}\nonumber\\
&&+\frac{1}{9}L_{n}^{2}L_{n-2+2(N-k)}^{3} 
+\frac{46}{125}L_{n-2}^{2}L_{n-1+2(N-k)}^{3}
+\frac{122}{225}L_{n-1}^{2}L_{n-1+2(N-k)}^{3}\nonumber\\
&&+\frac{29}{90}L_{n}^{2}L_{n-1+2(N-k)}^{3}+
\frac{59}{125}L_{n-1}^{2}L_{n+2(N-k)}^{3}
+\frac{6}{25}L_{n}^{2}L_{n+2(N-k)}^{3}\nonumber\\
&&+\frac{12}{125}L_{n}^{2}L_{n+1+2(N-k)}^{3}\nonumber\\
&&+\frac{12}{125}L_{n-3}^{3}L_{n-1+3(N-k)}^{2}+
\frac{6}{25}L_{n-2}^{3}L_{n-1+3(N-k)}^{2}
+\frac{29}{90}L_{n-1}^{3}L_{n-1+3(N-k)}^{2}\nonumber\\
&&+\frac{1}{9}L_{n}^{3}L_{n-1+3(N-k)}^{2}
+\frac{59}{125}L_{n-2}^{3}L_{n+3(N-k)}^{2}
+\frac{122}{225}L_{n-1}^{3}L_{n+3(N-k)}^{2}\nonumber\\
&&+\frac{8}{45}L_{n}^{3}L_{n+3(N-k)}^{2}
+\frac{46}{125}L_{n-1}^{3}L_{n+1+3(N-k)}^{2}
+\frac{8}{75}L_{n}^{3}L_{n+1+3(N-k)}^{2}\nonumber\\
&&+\frac{8}{125}L_{n}^{3}L_{n+2+3(N-k)}^{2}\nonumber\\
&&+\frac{24}{125}L_{n-3}^{4}L_{n+4(N-k)}^{1}
+\frac{393}{800}L_{n-2}^{4}L_{n+4(N-k)}^{1}
+\frac{17}{40}L_{n-1}^{4}L_{n+4(N-k)}^{1}\nonumber\\
&&+\frac{3}{32}L_{n}^{4}L_{n+4(N-k)}^{1} 
+\frac{58}{125}L_{n-2}^{4}L_{n+1+4(N-k)}^{1}
+\frac{71}{200}L_{n-1}^{4}L_{n+1+4(N-k)}^{1}\nonumber\\
&&+\frac{3}{40}L_{n}^{4}L_{n+1+4(N-k)}^{1}
+\frac{37}{125}L_{n-1}^{4}L_{n+2+4(N-k)}^{1}
+\frac{3}{50}L_{n}^{4}L_{n+2+4(N-k)}^{1}\nonumber\\
&&+\frac{6}{125}L_{n}^{4}L_{n+3+4(N-k)}^{1}\nonumber\\
&&+\frac{2}{25}L_{n-2}^{1}L_{n-2+N-k}^{1}L_{n-2+2(N-k)}^{3}
+\frac{1}{10}L_{n-1}^{1}L_{n-2+N-k}^{1}L_{n-2+2(N-k)}^{3}\nonumber\\
&&+\frac{1}{8}L_{n}^{1}L_{n-2+N-k}^{1}L_{n-2+2(N-k)}^{3}
+\frac{2}{15}L_{n-1}^{1}L_{n-1+N-k}^{1}L_{n-2+2(N-k)}^{3}\nonumber\\
&&+\frac{1}{6}L_{n}^{1}L_{n-1+N-k}^{1}L_{n-2+2(N-k)}^{3}
+\frac{2}{9}L_{n}^{1}L_{n+N-k}^{1}L_{n-2+2(N-k)}^{3}\nonumber\\
&&+\frac{11}{25}L_{n-1}^{1}L_{n-1+N-k}^{1}L_{n-1+2(N-k)}^{3}
+\frac{21}{40}L_{n}^{1}L_{n-1+N-k}^{1}L_{n-1+2(N-k)}^{3}\nonumber\\
&&+\frac{29}{45}L_{n}^{1}L_{n+N-k}^{1}L_{n-1+2(N-k)}^{3}
+\frac{12}{25}L_{n}^{1}L_{n+N-k}^{1}L_{n+2(N-k)}^{3}\nonumber\\
&&+\frac{6}{25}L_{n-2}^{1}L_{n-2+N-k}^{3}L_{n+4(N-k)}^{1}
+\frac{3}{10}L_{n-1}^{1}L_{n-2+N-k}^{3}L_{n+4(N-k)}^{1}\nonumber\\
&&+\frac{3}{8}L_{n}^{1}L_{n-2+N-k}^{3}L_{n+4(N-k)}^{1}
+\frac{23}{40}L_{n-1}^{1}L_{n-1+N-k}^{3}L_{n+4(N-k)}^{1}\nonumber\\
&&+\frac{2}{3}L_{n}^{1}L_{n-1+N-k}^{3}L_{n+4(N-k)}^{1}
+\frac{3}{8}L_{n}^{1}L_{n+N-k}^{3}L_{n+4(N-k)}^{1}\nonumber\\
&&+\frac{13}{25}L_{n-1}^{1}L_{n-1+N-k}^{3}L_{n+1+4(N-k)}^{1}
+\frac{23}{40}L_{n}^{1}L_{n-1+N-k}^{3}L_{n+1+4(N-k)}^{1}\nonumber\\
&&+\frac{3}{10}L_{n}^{1}L_{n+N-k}^{3}L_{n+1+4(N-k)}^{1}
+\frac{6}{25}L_{n}^{1}L_{n+N-k}^{3}L_{n+2+4(N-k)}^{1}\nonumber\\
&&+\frac{12}{25}L_{n-2}^{3}L_{n+3(N-k)}^{1}L_{n+4(N-k)}^{1}
+\frac{29}{45}L_{n-1}^{3}L_{n+3(N-k)}^{1}L_{n+4(N-k)}^{1}\nonumber\\
&&+\frac{2}{9}L_{n}^{3}L_{n+3(N-k)}^{1}L_{n+4(N-k)}^{1}
+\frac{21}{40}L_{n-1}^{3}L_{n+1+3(N-k)}^{1}L_{n+4(N-k)}^{1}\nonumber\\
&&+\frac{1}{6}L_{n}^{3}L_{n+1+3(N-k)}^{1}L_{n+4(N-k)}^{1}
+\frac{1}{8}L_{n}^{3}L_{n+2+3(N-k)}^{1}L_{n+4(N-k)}^{1}\nonumber\\
&&+\frac{11}{25}L_{n-1}^{3}L_{n+1+3(N-k)}^{1}L_{n+1+4(N-k)}^{1}
+\frac{2}{15}L_{n}^{3}L_{n+1+3(N-k)}^{1}L_{n+1+4(N-k)}^{1}\nonumber\\
&&+\frac{1}{10}L_{n}^{3}L_{n+2+3(N-k)}^{1}L_{n+1+4(N-k)}^{1}
+\frac{2}{25}L_{n}^{3}L_{n+2+3(N-k)}^{1}L_{n+2+4(N-k)}^{1}\nonumber\\
&&+\frac{3}{25}L_{n-2}^{1}L_{n-2+N-k}^{2}L_{n-1+3(N-k)}^{2}
+\frac{3}{20}L_{n-1}^{1}L_{n-2+N-k}^{2}L_{n-1+3(N-k)}^{2}\nonumber\\
&&+\frac{3}{16}L_{n}^{1}L_{n-2+N-k}^{2}L_{n-1+3(N-k)}^{2}
+\frac{3}{10}L_{n-1}^{1}L_{n-1+N-k}^{2}L_{n-1+3(N-k)}^{2}\nonumber\\
&&+\frac{3}{8}L_{n}^{1}L_{n-1+N-k}^{2}L_{n-1+3(N-k)}^{2}
+\frac{1}{3}L_{n}^{1}L_{n+N-k}^{2}L_{n-1+3(N-k)}^{2}\nonumber\\
&&+\frac{14}{25}L_{n-1}^{1}L_{n-1+N-k}^{2}L_{n+3(N-k)}^{2}
+\frac{53}{80}L_{n}^{1}L_{n-1+N-k}^{2}L_{n+3(N-k)}^{2}\nonumber\\
&&+\frac{8}{15}L_{n}^{1}L_{n+N-k}^{2}L_{n+3(N-k)}^{2}
+\frac{8}{25}L_{n}^{1}L_{n+N-k}^{2}L_{n+1+3(N-k)}^{2}\nonumber\\
&&+\frac{4}{25}L_{n-2}^{2}L_{n-1+2(N-k)}^{1}L_{n-1+3(N-k)}^{2}
+\frac{4}{15}L_{n-1}^{2}L_{n-1+2(N-k)}^{1}L_{n-1+3(N-k)}^{2}\nonumber\\
&&+\frac{1}{6}L_{n}^{2}L_{n-1+2(N-k)}^{1}L_{n-1+3(N-k)}^{2}
+\frac{2}{5}L_{n-1}^{2}L_{n+2(N-k)}^{1}L_{n-1+3(N-k)}^{2}\nonumber\\
&&+\frac{1}{4}L_{n}^{2}L_{n+2(N-k)}^{1}L_{n-1+3(N-k)}^{2}
+\frac{1}{6}L_{n}^{2}L_{n+1+2(N-k)}^{1}L_{n-1+3(N-k)}^{2}\nonumber\\
&&+\frac{17}{25}L_{n-1}^{2}L_{n+2(N-k)}^{1}L_{n+3(N-k)}^{2}
+\frac{2}{5}L_{n}^{2}L_{n+2(N-k)}^{1}L_{n+3(N-k)}^{2}\nonumber\\
&&+\frac{4}{15}L_{n}^{2}L_{n+1+2(N-k)}^{1}L_{n+3(N-k)}^{2}
+\frac{4}{25}L_{n}^{2}L_{n+1+2(N-k)}^{1}L_{n+1+3(N-k)}^{2}\nonumber\\
&&+\frac{8}{25}L_{n-2}^{2}L_{n-1+2(N-k)}^{2}L_{n+4(N-k)}^{1}
+\frac{8}{15}L_{n-1}^{2}L_{n-1+2(N-k)}^{2}L_{n+4(N-k)}^{1}\nonumber\\
&&+\frac{1}{3}L_{n}^{2}L_{n-1+2(N-k)}^{2}L_{n+4(N-k)}^{1}
+\frac{53}{80}L_{n-1}^{2}L_{n+2(N-k)}^{2}L_{n+4(N-k)}^{1}\nonumber\\
&&+\frac{3}{8}L_{n}^{2}L_{n+2(N-k)}^{2}L_{n+4(N-k)}^{1}
+\frac{3}{16}L_{n}^{2}L_{n+1+2(N-k)}^{2}L_{n+4(N-k)}^{1}\nonumber\\
&&+\frac{14}{25}L_{n-1}^{2}L_{n+2(N-k)}^{2}L_{n+1+4(N-k)}^{1}
+\frac{3}{10}L_{n}^{2}L_{n+2(N-k)}^{2}L_{n+1+4(N-k)}^{1}\nonumber\\
&&+\frac{3}{20}L_{n}^{2}L_{n+1+2(N-k)}^{2}L_{n+1+4(N-k)}^{1}
+\frac{3}{25}L_{n}^{2}L_{n+1+2(N-k)}^{2}L_{n+2+4(N-k)}^{1}\nonumber\\
&&+\frac{1}{5}L_{n-1}^{1}L_{n-1+N-k}^{1}L_{n-1+2(N-k)}^{1}L_{n-1+3(N-k)}^{2}
+\frac{1}{4}L_{n}^{1}L_{n-1+N-k}^{1}L_{n-1+2(N-k)}^{1}
L_{n-1+3(N-k)}^{2}\nonumber\\
&&+\frac{1}{3}L_{n}^{1}L_{n+N-k}^{1}L_{n-1+2(N-k)}^{1}L_{n-1+3(N-k)}^{2}
+\frac{1}{2} L_{n}^{1}L_{n+N-k}^{1}L_{n+2(N-k)}^{1}
L_{n-1+3(N-k)}^{2}\nonumber\\
&&+\frac{4}{5}L_{n}^{1}L_{n+N-k}^{1}L_{n+2(N-k)}^{1}L_{n+3(N-k)}^{2}
+\frac{2}{5}  L_{n-1}^{1}L_{n-1+N-k}^{1}L_{n-1+2(N-k)}^{2}
L_{n+4(N-k)}^{1}\nonumber\\
&&+\frac{1}{2}L_{n}^{1}L_{n-1+N-k}^{1}L_{n-1+2(N-k)}^{2}L_{n+4(N-k)}^{1}
+\frac{2}{3} L_{n}^{1}L_{n+N-k}^{1}L_{n-1+2(N-k)}^{2}
L_{n+4(N-k)}^{1}\nonumber\\
&&+\frac{3}{4} L_{n}^{1}L_{n+N-k}^{1}L_{n+2(N-k)}^{2}L_{n+4(N-k)}^{1}
+\frac{3}{5}  L_{n}^{1}L_{n+N-k}^{1}L_{n+2(N-k)}^{2}
L_{n+1+4(N-k)}^{1}\nonumber\\
&&+\frac{3}{5}L_{n-1}^{1}L_{n-1+N-k}^{2}L_{n+3(N-k)}^{1}L_{n+4(N-k)}^{1}
+\frac{3}{4}L_{n}^{1}L_{n-1+N-k}^{2}L_{n+3(N-k)}^{1}L_{n+4(N-k)}^{1}\nonumber\\
&&+\frac{2}{3}L_{n}^{1}L_{n+N-k}^{2}L_{n+3(N-k)}^{1}L_{n+4(N-k)}^{1}
+\frac{1}{2}L_{n}^{1}L_{n+N-k}^{2}L_{n+1+3(N-k)}^{1}L_{n+4(N-k)}^{1}\nonumber\\
&&+\frac{2}{5}L_{n}^{1}L_{n+N-k}^{2}L_{n+1+3(N-k)}^{1}L_{n+1+4(N-k)}^{1}
+\frac{4}{5}L_{n-1}^{2}L_{n+2(N-k)}^{1}L_{n+3(N-k)}^{1}
L_{n+4(N-k)}^{1}\nonumber\\
&&+\frac{1}{2}L_{n}^{2}L_{n+2(N-k)}^{1}L_{n+3(N-k)}^{1}L_{n+4(N-k)}^{1}
+\frac{1}{3}L_{n}^{2}L_{n+1+2(N-k)}^{1}L_{n+3(N-k)}^{1}
L_{n+4(N-k)}^{1}\nonumber\\
&&+\frac{1}{4}L_{n}^{2}L_{n+1+2(N-k)}^{1}L_{n+1+3(N-k)}^{1}L_{n+4(N-k)}^{1}
+\frac{1}{5}L_{n}^{2}L_{n+1+2(N-k)}^{1}L_{n+1+3(N-k)}^{1}
L_{n+1+4(N-k)}^{1}\nonumber\\
&&+L_{n}^{1}L_{n+N-k}^{1}L_{n+2(N-k)}^{1}
L_{n+3(N-k)}^{1}L_{n+4(N-k)}^{1}
\label{quin}    
\end{eqnarray}

\end{document}